\documentclass[11pt]{article}
\usepackage{amsmath,amssymb,epsf}
\usepackage[french]{babel}
\usepackage[applemac]{inputenc}
\advance \textwidth 2.5cm
\advance \textheight 3cm
\usepackage{amsfonts}
\usepackage{latexsym}
\newtheorem{theorem}{Th\'eor\`eme}

\newtheorem{corollaire}{Corollaire}

\newtheorem{prop}{Proprit}
\newtheorem{remarque}{Remarque}
{\unskip\nobreak\hfil\penalty50\hskip2em\null\nobreak\hfil%
$\Box$\parfillskip0pt\par\medskip}

\newcommand{\tore}{\mathbb T}
\newcommand{\Tr}{\mathrm {Tr}}

\title{Comportement asymptotique des polynmes orthogonaux associ\'es \`a un poids ayant un z\'ero d'ordre fractionnaire sur le cercle. Applications 
aux valeurs propres d'une classe de matrices al\'eatoires unitaires.}

\author{ Philippe Rambour\thanks{Universit\'{e} de Paris Sud,
      B\^atiment 425; F-91405
Orsay Cedex;
tel : 01 69 15 57 28 ; fax 01 69 15 60 19
      \mbox{e-mail : philippe.rambour@math.u-psud.fr}}
       \and Abdellatif Seghier\thanks{Universit\'{e} de Paris Sud,
        B\^atiment 425; F-91405
Orsay Cedex;
tel : 01 69 15 57 29 ; fax 01 69 15 72 34
       \mbox{ e-mail : abdelatif.seghier@math.u-psud.fr}}}

\begin{document}
\maketitle
  \renewcommand{\abstractname}{Rsum}
          \begin{abstract}
\textbf{ Comportement asymptotique des polynmes orthogonaux 
associ\'es \`a un poids ayant un z\'ero fractionnaire sur le cercle. Applications 
aux valeurs propres de certaines matrices al\'eatoires unitaires}\\
     Cet article s'int\'eresse aux polynmes orthogonaux 
     $\Phi _{N}$ correspondant aux poids
     de type $ (1-\cos \theta )^\alpha c$ o\`u $c$ est une fonction suffisamment r\'eguli\`ere et $\alpha \in ]-\frac{1}{2}, \frac{1}{2}[$. On d\'ecrit le comportement asymptotique des coefficients et de la valeur en $1$ de ces polynmes 
     et de leur d\'eriv\'ees lorsque $N$ tend vers l'infini. 
     Ceci nous permet d'obtenir une asymptotique du noyau de Christofel-Darboux 
     associ\'e \`a un tel poids et de calculer la loi conjointe des valeurs propres 
     d'une famille de matrices al\'eatoires unitaires. Les d\'emonstrations des r\'esultats 
     relatifs aux polyn\^omes orthogonaux sont essentiellement bas\'ees sur les propri\'et\'es de l'inverse de la matrice de Toeplitz de symbole $f$. 
        \end{abstract}
     \renewcommand{\abstractname}{Abstract}
          \begin{abstract}
\textbf{Asymptotic behavior of orthogonal polynomials on the circle, 
with respect  to a weight having a fractional zero on the torus. Applications to the eigenvalues of certain unitary random matrices. }\\
This paper is devoted to the orthogonal polynomial on the circle, with respect to
a weight of type $ f=(1-\cos \theta )^\alpha c$ where $c$ is a sufficiently smooth function and $\alpha \in ]-\frac{1}{2}, \frac{1}{2}[$.  We obtain an asymptotic expansion
of the coefficients of this polynomial and of $\Phi^{(p)}_{N}(1)$ for all integer $p$.
These results  allow us to obtain an asymptotic expansion of the associated Christofel-Darboux kernel, and to compute 
the distribution of the eigenvalues of a family of random unitary matrices.
 The proof of the resuts related with the orthogonal polynomials are essentialy based on the inversion of Toeplitz matrice associated to the symbol $f$. 

  \end{abstract}


\textbf{\large{Mots clef}}

\textbf{ Inversion des matrices de Toeplitz , matrices al\'eatoires unitaires, 
Polyn\^omes orthogonaux.}

\section{Introduction.}
Il est bien connu qu'une matrice alatoire est caractrise par la loi de probabilit de la 
distribution de ses valeurs propres. Lorsqu'il s'agit d'\'etudier des matrices alatoires unitaires, le cas le plus tudi est celui de l'ensemble circulaire de Dyson o cette densit est de la forme,
 pour une matrice unitaire $N\times N$ (voir \cite{Dy62a}, \cite{Dy62b}, \cite{Dy62c},
 \cite{TW00}) :
 $$ P_N (\theta_1, \cdots \theta_N) = c_N \prod_{1\le j<k \le N } \vert e^{i \theta_j} - e^{i\theta_k} \vert ^2.$$
 Il faut remarquer que dans le cas des matrices alatoires unitaires c'est bien s\^ur la
 distibution des arguments des valeurs propres qu'on tudie. \\
 Cet ensemble est considr comme une bonne reprsentation des phnomnes physiques qui ont motiv initialement l'tude des matrices alatoires. Il est  naturel de chercher \`a l'tendre  un ensemble plus vaste o la loi des valeurs propres est de la forme ( voir \cite{Nagao2} et \cite{TW01}) :
 \begin{equation} \label{Dyson}
  P_N (\theta_1, \cdots \theta_N) = \prod _{1\le j \le N}  f (\theta_j)
  \prod_{1\le j<k \le N } \vert e^{i \theta_j} - e^{i\theta_k} \vert ^2, 
  \end{equation} 
 (ensemble circulaire de Dyson gnralis) o\`u $f$ est une fonction positive, int\'egrable sur le tore. Alors que, dans le cas de l'ensemble circulaire de Dyson, la fonction de corrlation  des arguments des valeurs propres se calcule en fonction des polyn\^omes orthogonaux associs  la mesure de Haar du tore 
 (\cite{TW02}), dans le cas de l'ensemble circulaire de Dyson gnralis cette fonction de corrlation se calcule au moyen des 
 noyaux de Christoffel-Darboux $K_N$ associs aux polyn\^omes orthogonaux unitaires
$(\Phi_n)$ o\`u $1\le n\le N$ )associs au poids $f$. On d\'efinit g\'en\'eralement $K_{N}$ par  (voir \cite {SZEG}) :
$$ K_N (e^{i\theta}, e^{i\theta'}) =  \sqrt{f(e^{i \theta}) }
\sqrt {f(e^{i \theta'})}
\sum_{m=0}^{N-1} \frac{1}{h_m} \overline{ \Phi_m (e^{i \theta}) }
\Phi_m (e^{i \theta'}) 
$$ 
avec  
$$
\int_{- \pi} ^{\pi} f(e^{i \theta}) \overline{ \Phi_m (e^{i \theta}) } \Phi_m (e^{i \theta}) d\theta 
= h_m.
$$ 
On sait que, pour $z \not=0$, on peut \'ecrire aussi, si $\theta\not= \theta'$, et avec
$ \Phi_{N}^{*}(z) =z^N \bar \Phi _{N} \left( \frac{1}{z}\right),
 $
$$  K_N (e^{i\theta}, e^{i\theta'})=  \frac{1}{h_{N}} \sqrt{f(e^{i \theta}) }
\sqrt {f(e^{i \theta'})}
\frac {\overline{\Phi_{N} ^* (e^{i\theta})} \Phi_{N} ^* (e^{i\theta'})
-\overline {\Phi_{N}  (e^{i\theta}) } \Phi_{N} (e^{i\theta'})}
{ (1-e^{i(\theta'-\theta)} )}. $$
Dans le cas unitaire, le lien entre ce noyau et la fonction de corr\'elation est rappel\'e dans le cours de la d\'emonstration.
La fonction de corr\'elation \'etant d\'efinie, un des probl\`emes relatifs aux matrices al\'eatoires
est de d\'eterminer la distribution des valeurs propres lorsque l'entier $N$ tend vers l'infini. Ce probl\`eme se subdivise 
en fait en deux. D'une part il s'agit de d\'eterminer la distribution des arguments de ces valeurs propres 
sur un intervalle $]\epsilon, 2\pi -\epsilon[$ avec $\epsilon>0$  pour les matrices al\'eatoires unitaires ou la 
distribution des valeurs propres elles-m\^emes dans 
l'intervalle $]0,+ \infty[$ pour les matrices al\'eatoires hermitiennes  (voir \cite{TW94a}, \cite {TW01}, \cite{TW00}).
D'autre part on recherche cette distribution 
au voisinage de z\'ero pour les matrices al\'eatoires unitaires, ou au voisinage de z\'ero ou de plus l'infini 
pour les matrices al\'eatoires hermitiennes. \\
Dans le cas de matrices  al\'eatoires Hermitiennes  Gaussiennes et gr\^ace \`a des rsultats   
qui  ont t obtenus sur les polyn\^omes orthogonaux qui leur sont associ\'es (on peut consulter, par exemple, \cite {BI99}, \cite {D01}, \cite {Ku04}, \cite{De98} \ldots), on dispose de nombreux renseignements sur la loi des valeurs propres `` au bord'' 
du spectre (c'est \`a dire les valeurs propres proches de z\'ero ou tendant vers l'infini).\\
Dans ce cadre Kuijlaars et Vanlessen  \cite {Ku04} tudient des matrices alatoires Gaussiennes 
(hermitiennes) dont la loi des valeurs propres est donne par 
$$ P^{(N)} (x_{1},x_{2},\cdots x_{N} ) = \frac{1}{Z_{n}} \prod_{j=1}^N
w_{N} (x_{j} )\prod _{i<j} \vert x_{i} -x_{j }\vert ^2,$$
o\`u  $w_{N}$ un poids de la forme 
$ w_{n}(x) = \vert x \vert ^{2\alpha} e^{-NV(x)}, \quad \mathrm{pour} \, x\in \mathbb R $ et  $V$ une fonction potentielle qui peut par exemple 
\^etre une fonction polyn\^ome.
Les auteurs de l'article montrent que la fonction de corrlation des valeurs propres autour de l'origine est donn\'ee par 
un noyau de la forme 
$$ \mathbb J _{\alpha}^\alpha (u,v) = \pi \sqrt u \sqrt v \frac{J_{\alpha+1/2}(\pi u) J_{\alpha-1/2}(\pi v) -J_{\alpha-1/2}(\pi u) J_{\alpha+1/2}(\pi v) }
{2(u-v)},$$
o\`u les fonctions $J_{\alpha+1/2}$ et  $J_{\alpha-1/2}$ sont les fonctions de Bessel usuelles.\\
Dans le prsent  article, nous obtenons un r\'esultat du m\^eme type relativement \`a une classe de matrices al\'eatoires unitaires. En fait nous \'etudions les coefficients des polynmes orthogonaux $\Phi_{N}$ de degr\'e $N$ associ\'es au symbole
\begin{equation} \label{ici}
 f(e^{i \theta}) = \vert 1- e^{i\theta} \vert ^{2\alpha } c(e^{i \theta})
 \end{equation}
avec $c$ une fonction rgulire sur le tore et $ \alpha \in ] - \frac{1}{2}, \frac{1}{2} [$. En effet nous obtenons dans un premier temps une expression asymptotique quand $N$ tens vers l'infini des coefficients des polyn\^omes orthogonaux puis nous en  
d\'eduisons 
 la loi des arguments des valeurs propres au voisinage de l'origine quand $N$ tend vers l'infini.\\
Nous donnons d'abord (thormes \ref {Theo1} et \ref{Theo1bis}) 
l'asymptotique quand $N$ tend vers l'infini des coefficients des polyn\^omes  $\Phi^*_{N}$ et $\Phi_{N}$. 
Pour ce faire nous utilisons des mthodes bases sur des rsultats antrieurs (se rapporter  \cite{RS08} ou \cite{RS10}). 
   Si l'on pose $\displaystyle{\Phi_N (e^{i \theta}) = \sum_{u=0}^N \tilde\beta_u e^{iu \theta}}$ 
notre nonc se subdivise en trois parties qui correspondent  trois zones distinctes de l'ensemble des
indices. Nous donnons d'abord les indices qui correspondent "aux bords" de l'ensemble des indices,
c'est  dire des  entiers $u$ tels que 
$\displaystyle{\lim_{N\rightarrow + \infty} \frac{u}{N} =0}$ ou 
$\displaystyle{\lim_{N\rightarrow + \infty} \frac{u}{N} =1}$ et les indices qui se trouvent "au coeur"
de l'ensemble des indices, ce qui correspond aux entiers $u$ tels que 
$\displaystyle{\lim_{N\rightarrow + \infty} \frac{u}{N} =x}$ avec $0<x<1$.
Ces expressions nous permettent ensuite de calculer 
l'asymptotique  de $\Phi^{(p)}_N (1)$ et $\Phi^{*(p)}_N (1)$
(voir les thormes \ref{Theo1} et \ref{Theo22}), pour un entier $p$ fix.
Nous pouvons alors calculer une reprsentation asymptotique du noyau de 
Christoffel-Darboux d'ordre $N$ associ au symbole 
$\theta \rightarrow \vert 1- e^{i\theta} \vert ^{2\alpha } c(e^{i \theta})$.
Si $-\pi < u,v \le \pi$ et $u \not= v $, on a 
\begin{equation} \label{groscas}
   K_{N} (e^{i u/N},e^{iv/N} ) 
= \frac{ N (uv)^\alpha}{ \Gamma^2 (\alpha) c_{1} (1)} \frac  {\rho (\alpha,-u) \rho(\alpha,v) 
- \rho (\alpha,u) \rho(\alpha,-v) 
 e^{i(v-u)}} {i(u-v)} +o(N)
 \end{equation}
o $\rho$ est une fonction qui sera prcise plus loin (voir les thormes \ref{proba1} et
\ref{proba2}), et dont l'expression est diff\'erente selon que $\alpha \in ]- \frac{1}{2},0[$ 
ou que $\alpha\in ]0, \frac{1}{2}[$. \\
 Notons maintenant $ \theta_1^N, \cdots \theta_N^N$ les arguments contenus 
dans l'intervalle $]-\pi, \pi]$ d'une matrice alatoire $U_N$ de taille $N \times N$ 
suivant la loi donne dans l'\'egalit\'e (\ref{Dyson}) avec un poids $W$ comme pr\'ecis\'e en (\ref{ici}). Si l'on pose 
$F_N = \displaystyle{\sum_{i=1}^N \delta _{\theta_i^N}}$ les rsultats prcdents permettent d'tablir, si $- \pi <u<v\le \pi $ (voir les thormes \ref{proba1} et
\ref{proba2}) 
\begin{itemize} 
\item [i)]
$$ \displaystyle{ \lim_{N \rightarrow + \infty}
 P \Bigl( F_N [ \frac{u}{N}, \frac{v}{N}] = m \Bigr) }= F_{\infty} (a,b,m)$$
o $F_{\infty} $ est une  fonction dfinie au moyen du dterminant d'un oprateur de 
Fredholm de noyau $K= \lim_{N\rightarrow + \infty} \frac{K_{N}} {N}$. 
\item[ii)]
Si $q$ est un entier suprieur ou gal  2 
$$  \displaystyle{ \lim_{N \rightarrow + \infty}
 P \Bigl( F_N [ \frac{u}{N^q}, \frac{v}{N^q}] = m \Bigr) }=0; $$
 \end {itemize}
 Heuristiquement, cela signifie que lorsque $N$ tend vers l'infini il existe des valeurs propres de la forme
$e^{i \theta_N}$ avec $\theta_N = O(\frac{1}{N})$, mais qu'il n'en existe pas de la forme 
$e^{i \theta'_N}$ avec $\theta'_N = O(\frac{1}{N^q})$ pour $q>1$.\\
La premire partie de nos rsultats se rapproche en le compl\'etant dutravail de Martinez-Finkelstein, Mac-Laughin et Saff
qui dans \cite{ML3} tudient une famille de polyn\^omes orthogonaux $\tilde \Phi_N $ orthogonaux pour le 
poids \begin{equation} \label{poids}
W(z) = w(z) \prod_{k=1}^m \vert z- a_k \vert ^{2 \beta_k} \quad z \in \tore,
\end{equation}
o\`u la fonction $w$ est suffisamment rgulire et o $W$ vrifie les conditions de Szeg\"o
\begin{equation} \label {C1}
\int_{\tore} W(z) dz >0
\quad 
\mathrm {et} \quad 
\int_{\tore}\ln ( W(z) ) dz > - \infty.
\end{equation}
D ans l'article  (\cite{ML3}) ces auteurs donnent notamment 
le comportement asymptotique de $\Phi_N (z)$ pour $z \not=a_k \quad 1\le k \le m$, et 
tudient galement les zros des polyn\^omes $\tilde \Phi_N$, leur principal outil de dmonstration tant la mthode de Riemann-Hilbert.\\
A la fin de cet article, nous proposons un appendice dans lequel nous donnons des 
tableaux de valeurs numriques permettant de comparer 
$\Phi_N (1)$ et les quantits donnes dans les thormes  
\ref{Theo2} et \ref{Theo22} pour des valeurs de $\alpha$ \'egales \`a 
$-0,275$, $-0,150$, $0,025$, $0,1$, $0,225$.

\section{Notations et d\'efinitions. Enonc\'e des r\'esultats}
Avant d'\'enoncer nos principaux r\'esultats, nous introduisons
des notations et d\'efinitions.
\subsection{Notations. D\'efinitions}
Etant donn\'ee  $f \in L^1(\mathbb{T})$ une fonction d\'efinie sur
le cercle unit\'e, on appelle  matrice de Toeplitz de taille
$(N+1)\times (N+1)$ associ\'ee au symbole $f$ la  matrice
\begin{equation}
T_N(f)=\left(\widehat{f}(i-j\right)_{0\le i,j\le N}
\end{equation}
o\`u
\begin{equation}
\widehat{f}(k)=\int_{-\pi}^{\pi}
f(e^{i\theta})e^{-ik\theta} \frac{d\theta}{ 2\pi},~k\in\mathbb{Z}
\end{equation}
d\'esigne la suite des  coefficients de Fourier de $f$. \\
 Le lien entre les matrices de Toeplitz et les
polyn\^omes orthogonaux est bien connu
(voir, par exemple, \cite{Ld}).
Ainsi nous avons
\begin{equation}
\Phi^*_N(z)=\sum_{k=0}^N\frac {(T_N(f))^{-1}_{k+1,1}}
{(T_N(f))^{-1}_{1,1}}z^{k},~\mid z\mid =1.
\end{equation}
Les polyn\^omes $\Phi_n ^* \sqrt{(T_n(f))^{-1}_{1,1}}$ sont souvent appels polyn\^omes prdicteurs.  Comme le montre la relation ci dessus, leurs coefficients sont, \`a normalisation pr\`es, ceux de la premire colonne de la matrice ${T_N (f)}^{-1}.$ Ce sont d'abord ces coefficients que nous donnons dans les thormes \ref{Theo1} et \ref{Theo1bis}.\\
D'autre part la relation 
  \begin{equation} \label{predi}
  \Phi_N ^* (z) = z^N \bar \Phi_N (\frac{1}{z}),
\end{equation}
implique que les coefficients des polyn\^omes orthogonaux correspondent,  une normalisation prs, aux termes de la dernire colonne de la matrice $T_N^{-1}(f)$. Nous utiliserons la notation, pour tout $k$ dans 
$\mathbb Z$
$$
\chi^k : e^{i\theta}\mapsto \chi^k(\theta)=e^{ik\theta}.
$$
On notera $H_+$ (resp.$H_-$) le sous-espace ferm\'e de
$L^2(\mathbb{T})$ constitu\'e des fonctions $u$ dont les
coefficients de Fourier $\widehat{f}(k)=0$ quand $k<0$ (resp.
$k\ge 0)$. On notera par $\pi_\pm : L^2(\mathbb{T})\mapsto H_\pm$
les projecteurs orthogonaux respectifs.
\par
\noindent Nous rappelons  la d\'efinition de l'espace de Beurling.
On commence par introduire la notion de poids de Beurling : une
suite $(\mu_j)_{j\in\mathbb{Z}}$ est dite poids de Beurling si
elle satisfait aux trois propri\'et\'es suivantes
\begin{equation}{}
\left\{
\begin{array}{lll}
 (i) & \mu_j\ge 1,& j \in \mathbb{Z}\\
 (ii) & \mu_j=\mu_{-j},& j\in \mathbb{Z}\\
(iii) & \mu_{j+k}\le \mu_j\mu_k,&j,k\in\mathbb{Z}.
\end{array}
\right.
\end{equation}
La classe de Beurling ${\cal W}_\mu$ associ\'ee au poids
$\mu=(\mu_j)$ est alors d\'efinie comme suit
\begin{equation}
{\cal W}_\mu=\lbrace w \in
L^1(\mathbb{T})~:~\sum_{j\in\mathbb{Z}}\mu_j \mid
\widehat{w}(j)\mid <\infty\rbrace.
\end{equation}
Soit $f$ un symbole tel que 
$$f(e^{i \theta})=\mid 1-e^{i \theta}\mid^{2\alpha}c(e^{i \theta})$$ 
o\`u $c$ est une fonction suppos\'ee \^etre strictement
positive sur cercle unit\'e et appartenant \`a un espace de
Beurling. Nous noterons $g_\alpha \in H_+$ la fonction analytique
dans le disque unit\'e ouvert, dite fonction ext\'erieure,
v\'erifiant
$$
f=\mid g_\alpha \mid^2,
$$
et $c_1$ la fonction de $H^+$ vrifiant 
$$ c = c_1 \bar c_1$$
nous noterons par $(\beta_k^{(\alpha)})$ la suite des coefficients
de Fourier de $1/g_\alpha$,
et enfin pour allger les notations nous supposerons 
$(\beta_k^{(0)}) =1$.

\subsection{Enonc\'e des r\'esultats.} 
\subsubsection {Polyn\^omes orthogonaux}
Nous nous int\'eressons, en premier lieu, au comportement
asymptotique de la premi\`ere colonne de l'inverse d'une matrice
de Toeplitz $T_N(f)$. 
Nous donnons d' abord un nonc qui porte sur les coefficients situs sur les " bords "
de la premire colonne de l'inverse (voir \cite{KaRS}).
\begin{theorem} \label{Theo1}
Soit $T_N(f)$ une
matrice de Toeplitz dont le symbole $f$ s'\'ecrit
$$
f =\mid 1-\chi\mid^{2\alpha} c
$$
o\`u $\alpha \in \mathbb{R},~\mid \alpha \mid < \frac {1}{ 2}$ et o\`u
la fonction  $c>0$  est suppos\'ee appartenir \`a une classe de
Beurling $\mathcal {W}_{\mu}$ avec $\mu_{j}\ge \frac{3}{2}$ pour tout entier $j$.  Pour tout entier naturel $k$, tel que 
$\frac {k}{ N}\rightarrow 0$
 lorsque $N\rightarrow \infty$ on a alors
\par
\noindent i)
\begin{equation}
T_N(f)^{-1}_{k+1,1}=\left(\beta^{(\alpha)}_k- \frac {\alpha^2}{ N}\beta^{(\alpha+1)}_k
\right)+ R_N,~~ N\rightarrow\infty
\end{equation}
avec $R_N = O (\frac{k^{\alpha+1}}{N^2} )$ si $0<\alpha <\frac{1}{2},$
et avec $ R_N = O (\frac{k}{N^2} ) $ siÊ$-\frac{1}{2}<\alpha <0.$
\par
\noindent ii)
\begin{equation}
\left(T_{N}^{-1}(f)\right)_{N+1-k,1}=\left( 
\frac {c_ (1))}
{\bar c_1(1)} \beta_k^{(\alpha+1)}\right) \frac {\alpha}
{N}+o( \frac {1}{ N}),~N\rightarrow \infty.
\end{equation}
La 
convergence est uniforme sur tout intervalle $[0, [N \epsilon]]$ et sur $(N, N -[N \epsilon])$
pour $\epsilon$ assez petit.
\end{theorem}
On a d'autre part l'nonc suivant qui porte sur les lments du "coeur" de la premire colonne de l'inverse (\cite{RS08} et \cite{RS10})
\begin {theorem} \label{Theo1bis}
Si $0<x<1$ est un r\'eel donn\'e, on a, avec les m\^emes hypothses que pour le thorme
\ref{Theo1} 
\begin{equation}
c_{1}(1) \left(T_{N}^{-1}(f)\right)_{[Nx]+1,1}= 
K_\alpha(x) N^{\alpha-1}+o(N^{\alpha-1}),~N\rightarrow\infty
\end{equation}
o\`u
$$
K_\alpha(x)= \frac {1}{ \Gamma(\alpha)}x^{\alpha-1}(1-x)^\alpha.
$$
 Le thorme \ref{Theo1bis} la convergence est uniforme sur tout intervalle 
$[\delta _{1}, \delta _{2}]$ pour $0<\delta _{1} <\delta_{2}<1$.
\end{theorem}
 Avec les hypothses que nous nous sommes donnes, les termes calculs ci-dessus correspondent aux coefficients de
$\Phi^{*} _N$. Ceux de $\Phi_N$ en sont alors dduits  partir de la relation 
\ref{predi}. Ce sont ces coefficients que nous donnons maintenant. 

\begin{corollaire} \label{uno}
Posons  $\Phi_N (z) = \displaystyle{ \sum_{u=0}^{N+1} \tilde \beta_u z^u}$ le polynme orthogonal 
associ au symbole $f =\vert 1- \chi \vert ^{2 \alpha} c$. Avec les mmes hypothses que pour les thormes 
\ref{Theo1} et \ref{Theo1bis} on a, 
si $k/N \rightarrow 0$ quand $N \rightarrow + \infty$  
\begin {enumerate}
\item
$$\tilde\beta_{N-k}=\left(\overline{\beta_k^{(\alpha)}}- \frac {\alpha^2}{ N} \overline{\beta_k^{(\alpha+1)}}
\right)+ R_N,~~ N\rightarrow\infty
$$
avec $R_N = O (\frac{k^{\alpha+1}}{N^2} )$ si $0<\alpha <\frac{1}{2},$
et avec $ R_N = O (\frac{k}{N^2} ) $ siÊ$-\frac{1}{2}<\alpha <0.$
\item
$$\tilde \beta_k =\left( 
 \frac{\bar c_ 1(1)}
{ c_1 (1)} \overline{\beta_{N-k}^{(\alpha+1)}}\right) \frac{\alpha}
{N}+o(\frac{1}{ N}),~N\rightarrow \infty .$$
\end{enumerate}
De plus, si $0<x<1$ est un r\'eel donn\'e, on a alors
\begin{equation}
\tilde \beta_{[Nx]}= K_\alpha(1-x) N^{\alpha-1}+o(N^{\alpha-1}),~N\rightarrow\infty
\end{equation}
\end{corollaire}
Comme nous l'avons rappel dans l'introduction, le calcul du noyau de Christoffel-Darboux 
que nous nous proposons d'effectuer pour obtenir les noyaux annoncs nous 
oblige  dterminer avec prcision les quantits $\Phi_N^{(p)}(1)$ et 
$\left(\Phi_N ^*\right)^{(p)} (1)$ pour tout entier naturel $p$.
C'est ce que nous faisons dans les thormes \ref{Theo2} et \ref{Theo22}. 
Il est  noter que ces rsultats, bien que donns seulement dans le cas d'une seule singularit, prolongent les rsultats de \cite{ML3}, puisque dans cet article les valeurs prises en les singularits par les polyn\^omes orthogonaux sont inaccessibles.
\begin{theorem} \label{Theo2}
Avec les hypoth\`eses des th\'eor\`emes (\ref{Theo1})  et (\ref{Theo1bis}) on obtient
\begin{itemize}
 \item 
  Si ${-1\over 2} < \alpha <0 $ on a
 \begin{equation}
\Phi_N^*(1)= 
 {N^\alpha \over \Gamma(\alpha)}\frac {1}{c_1(1)} 
\left(\int_0^1 x^{\alpha -1} \left( (1-x)^\alpha -1\right) dx + \frac{1}{\alpha}
\right) +o(N^{\alpha}).
 \end{equation}
 \item
  D'autre part, pour un entier $j$ strictement positif fix\'e, dans le cas   ${-1\over 2} < \alpha <0 $ 
 et pour un entier $j$ positif fix\'e dans le cas $0<\alpha<{1\over 2}$ on a
 \begin{equation}
\left(\Phi_N^*\right) ^{(j)}(1)=  {N^{\alpha+j} \over \Gamma(\alpha)} {1\over
{c_1 (1)}} \int_0^1
x^{\alpha-1+j} (1-x)^\alpha  dx+o(N^{\alpha+j}).
 \end{equation}
  \end{itemize}
\end{theorem}
\begin{remarque}
On peut remarquer que pour un entier $j$ strictement positif et suffisamment petit devant $N$
toujours avec l'hypoth\`ese $-\frac{1}{2} <\alpha < \frac{1}{2}$ on a en fait 
$$ \left(\Phi_N^*\right) ^{(j)}(1) = N^{\alpha+j} \frac{\Gamma (\alpha+j)\Gamma(\alpha +1)}
{\Gamma (2 \alpha +j+1)} \frac{1} {\Gamma (\alpha) c_1(1)}+o(N^{\alpha+j}) .$$
\end{remarque}
On obtient le m\^eme genre d'\'enonc\'e pour les polyn\^omes orthogonaux. 
\begin{theorem} \label{Theo22}
Avec les hypoth\`eses des th\'eor\`emes \ref{Theo1}  et  \ref{Theo1bis} on obtient
\begin{itemize}
\item
 Si $-\frac{1}{ 2} < \alpha <0 $ et si $ j$ un entier positif 
on a  
 \begin{equation}
\Phi_N ^{(j)} (1)=
\frac {N^{\alpha+j} }{ \Gamma(\alpha)}\frac {1}
{\bar c_1(1)} 
\left(\int_0^1 x^{\alpha -1} \left( (1-x)^{\alpha+j}-1 \right) dx + \frac{1}{\alpha}
\right) +o(N^{\alpha}).
 \end{equation}
 \item
D'autre part, pour un entier $j$  positif fix\'e dans le cas $0<\alpha<{1\over 2}$ nous pouvons \'ecrire 
 \begin{equation}
\left(\Phi_N \right) ^{(j)}(1)={N^{\alpha+j} \over \Gamma(\alpha)} \frac{1}{\bar c_1(1)} \int_0^1
x^{\alpha-1} (1-x)^{\alpha +j} dx+o(N^{\alpha+j}).
 \end{equation}
 \end{itemize}
\end{theorem}
\begin{remarque}
On peut l\`a aussi remarquer que pour un entier $j$ strictement positif 
suffisamment petit devant $N$ et toujours avec $\frac{1}{2}>\alpha>0$ on a
$$ \left(\Phi_N \right) ^{(j)}(1) =N^{\alpha+j}  \frac{\Gamma (\alpha)\Gamma(\alpha +j+1)}
{\Gamma (2 \alpha +j+1)} \frac{1} {\Gamma (\alpha) \bar c_1(1)} + o(N^{\alpha+j} ).$$
\end{remarque}
\subsubsection{Matrices alatoires unitaires.}

Soit maintenant $U_{N}$ une matrice al\'eatoire unitaire de taille $N$.
 Notons $\theta_i^N$ $1 \le i \le N$ les arguments des valeurs propres de $U_{N}$ que l'on suppose class\'es dans l'ordre croissant dans $] - \pi , + \pi]$.
  On suppose que ces arguments ont pour densit\'e de probabilit\'e 
  $$ P_N (\theta_1, \theta_2, \cdots, \theta_N) = \displaystyle{\prod_{j=1}^N
 f(e^{i\theta_j}) \prod_{j<i} ^N 
\vert e^{i \theta_j} - e^{i \theta_i} \vert ^2}$$
avec $ f(e^{i \theta}) = \vert 1- e^{i \theta} \vert ^{2 \alpha} c_1(e^{i \theta})$ 
et o\`u $c_{1}$ est dans une classe de Beurling d'indice sup\'erieur 
ou \'egal \`a $\frac{3}{2}$.
Nous posons ensuite
$F_N = \displaystyle{\sum_{i=1}^N \delta _{\theta_i^N}}$ 
o\`u $\delta _{\theta_i^N}$ d\'esigne la mesure de Dirac de support $\theta_i^N$.
Dans la suite nous supposerons enfin que $u$ et $v$ sont deux rels tels que 
$-\pi<u\le v< \le \pi$. On a alors les \'enonc\'es suivants qui prcisent la loi au voisinage de $1$ des valeurs propres des matrices alatoires unitaires tudies
 dans cet article.
\begin{theorem} \label{proba1}
Soit $ 0<\alpha<\frac{1}{2}$. Posons, pour tout $u \in \mathbb R$ 
$$
\psi (\alpha,u) = \int_{0}^1 x^{\alpha-1} (1-x)^\alpha e^{iux} dx,
\quad  \tau(\alpha,u) =\int_{0}^1 x^{\alpha} (1-x)^\alpha e^{iux} dx
$$
 Nous pouvons alors \'ecrire, 
avec les hypoth\`eses ci-dessus, et pour tout entier positif $m$ 
\begin{equation} \label {chebello}
 \displaystyle{ \lim_{N \rightarrow + \infty}
 P \Bigl( F_N [ \frac{u}{N}, \frac{v}{N}] = m \Bigr) }
= \frac{(-1)^m} {m!} \left( \frac{d}{d\gamma}\right)^m 
\det [ (Id - \gamma \mathcal K)\vert _{L^2(I)}] \Big \vert _{\gamma =1}
\end{equation}
o\`u $I=[u,v]$et o $\mathcal K$ est l'op\'erateur de Fredhom de noyau $K $
avec 
\begin{itemize} 
\item 
si $\theta \not= \theta'$ 
$$ 
 (\theta,\theta') \rightarrow 
 \frac{\theta^\alpha \theta^{\prime\alpha}}{\Gamma^2(\alpha) c(1) e^{-i \theta}}
 \Bigl (\frac{ e^{i \theta}\psi (\alpha, -\theta)\psi (\alpha, \theta') - \psi (\alpha, \theta)\psi (\alpha,- \theta') 
 e^{i\theta'}} {i(\theta'-\theta)}\Bigr ) .
$$
\item
et si $\theta=\theta'$, 
$$ \theta \rightarrow 
\left( -2 \Re \left(\psi (\alpha,\theta) \tau(\alpha,-\theta)\right  )+
 \vert \psi (\alpha, \theta) \vert ^2\right) \frac{\theta^{2 \alpha}} {\Gamma^2(\alpha) c (1)}.
$$
\end{itemize}
\end{theorem}
\begin {remarque}
On a en particulier le r\'esultat
\begin{equation}
 \displaystyle{ \lim_{N \rightarrow + \infty}
 P \Bigl( F_N [ \frac{u}{N}, \frac{v}{N}] = 0 \Bigr) } = \det [ (Id-\mathcal K)\vert _{L^2(I)}] 
\end{equation}
\end{remarque}
Le th\'eor\`eme suivant permet d'envisager le cas o\`u $\alpha $ est n\'egatif.
\begin{theorem} \label{proba2}
Soit $ -\frac{1}{2}<\alpha<0$. Posons, pour tout $u,v \in \mathbb R$ 
$$
\tilde \psi (\alpha,u) = \int_0 ^{1} x^{\alpha-1}( (1-x) ^\alpha e^{-i u x}-1) dx
+ \frac{1}{\alpha}.$$
 Nous pouvons alors \'ecrire, 
avec les hypoth\`eses du th\'eor\`eme \ref{proba1}, 
\begin{equation} \label {chebello2}
 \displaystyle{ \lim_{N \rightarrow + \infty}
 P \Bigl( F_N [ \frac{u}{N}, \frac{v}{N}] = m \Bigr) }
= \frac{(-1)^m} {m!} \left( \frac{d}{d\gamma}\right)^m 
\det [ (Id - \gamma \tilde {\mathcal K})\vert _{L^2(I)}] \Big \vert _{\gamma =1}
\end{equation}
o\`u $I=[u,v]$ o $\tilde {\mathcal K}$ est l'op\'erateur de Fredhom de noyau $ \tilde K $
avec 
\begin{itemize}
\item
si $\theta\not=\theta'$ 
$$
 (\theta,\theta') \rightarrow 
\frac{\theta^\alpha \theta^{\prime\alpha}}{\Gamma^2(\alpha) c(1) e^{-i \theta}} 
\Bigl ( \frac{e^{i \theta}\tilde\psi (\alpha, -\theta)\tilde\psi (\alpha, \theta') -
 \tilde \psi (\alpha, \theta) \tilde \psi (\alpha,- \theta') 
 e^{i\theta'}} {i(\theta'-\theta)})\Bigr ).
$$ 
\item
Ou encore si $\theta=\theta'$, 
$$ \theta \rightarrow 
\left( -2 \Re \left(\tilde \psi (\alpha,\theta) \tau(\alpha,-\theta)\right  )+
 \vert \tilde \psi (\alpha, \theta) \vert ^2\right) \frac{\theta^{2 \alpha}}
  {\Gamma^2(\alpha) c (1)}.
$$
\end{itemize}
\end{theorem}
\begin{remarque}
On a en toujours le cas particulier int\'eressant suivant 
\begin{equation}
 \displaystyle{ \lim_{N \rightarrow + \infty}
 P \Bigl( F_N [ \frac{u}{N}, \frac{v}{N}] = 0 \Bigr) } = \det [ (Id- \tilde {\mathcal K})\vert _{L^2(I)}] 
\end{equation}
\end{remarque}
 \begin{theorem}\label{proba3}
 Pour tous entiers naturels $q$ et $m$ v\'erifiant
$q>1$ et $N>m\ge1$  et pour tout r\'eel $\alpha$ tel que
$\frac{-1}{2} <\alpha <\frac{1}{2}$ on a, toujours avec les m\^emes hypoth\`eses qu'au 
th\'eor\`eme \ref{proba1} 
\begin{equation} \label {chebello3}
 \displaystyle{ \lim_{N \rightarrow + \infty}
 P \Bigl( F_N [ \frac{u}{N^q}, \frac{v}{N^q}] = m \Bigr) }=
O\left(N^{(1-p)}\right).
\end{equation}
\end{theorem}

\section{Preuves des r\'esultats}
\subsection{Preuve des th\'eormes \ref{Theo1} et \ref{Theo1bis} }
Les points $i)$ et $ii)$ du thorme \ref{Theo1} ont \'et\'e d\'emontr\'es dans \cite{KaRS}.
La preuve est bas\'ee sur  une formule explicite de l'inverse que nous avons obtenue dans un prcdent travail  (\cite{RS08} ou \cite{RS10}). Nous pouvons en effet crire 
 pour tout entier $k$, $0\le k \le N+1$,
l'expression 
\begin{equation} \label{superformule}
(T_N(f))^{-1}_{k+1,1}= \bar\beta^{(\alpha)}_0\left(
\beta_k^{(\alpha)}-\frac{1}{N} \left( \sum_{u=0}^k
\beta_u^{(\alpha)}F_{N,\alpha}\left((k-u)/N\right) \right)\right) \left(1+o(1)\right)
\end{equation}
 avec 
 $$
F_{N,\alpha} (u/N) =\sum_{m=0}^\infty
F_{2m,N}^{(\alpha)}(u/N)\left( {\sin{(\pi \alpha)}\over \pi}
\right)^{2m+2} 
$$
o\`u on  a pos\'e, pour $z\in [0,1]$
\begin{equation}\label{Wu}
\begin{array}{lll}
 F_{2m,N}^{(\alpha)}(z)&=&\displaystyle{\sum_{j_0}^\infty{1\over j_0+N_\alpha}\sum_{j_1=0}^\infty{1\over
j_1+j_0+N_\alpha} \ldots+
\sum_{j_{2m-1}=0}^\infty{1\over j_{2m-1}+j_{2m-2}+N_\alpha}}+\\
\displaystyle&+&\sum_{j_{2m}=0}^\infty{1\over
j_{2m}+j_{2m-1}+N_\alpha}{1\over ( j_{2m}+ \alpha+1)/N - z}
\end{array}
\end{equation}
o\`u $N_\alpha=N+1+\alpha$. 
Dans (\cite{I2}), Inoue montre que
${1\over N} F_{2m,N}^{(\alpha)}(0)$ est approxim\'e par  une
int\'egrale dont la valeur est $\alpha^2$. En combinant ce
r\'esultat \`a un  th\'eor\`eme d'approximation de Bleher (\cite{BL})
 permettant d'estimer l'\'ecart entre
l'int\'egrale et la somme, on obtient le d\'eveloppement
\begin{equation}
\sum_{m=0}^\infty F_{2m,N}^{(\alpha)}(0)\left( {\sin{(\pi
\alpha)}\over \pi} \right)^{2m+2}   ={\alpha^2 \over N}+O({1\over
N^2}).
\end{equation}
Ces rsultats permettent alors d'obtenir le point i) du thorme. Le point ii)
est une consquence de la relation de r\'ecurrence de Szeg\"o (voir \cite{GS}) :
 pour $\mid z\mid=1$ 
\begin{equation}
\left\{
\begin{array}{lll}
\Phi_{N+1}(z)&=&z\Phi_N(z)-\overline{\gamma_N}\Phi_n^*(z)\\
\Phi_{N+1}^*(z)&=&\Phi_n^*(z)-\gamma_N z\Phi_n(z),
\end{array}
\right.
\end{equation}
o\`u $\bar \gamma_{N} = - \overline{ \Phi_{n+1}(0)} $ est appel\'e coefficient de Verblunsky. 
L'uniformit est une consquence de l'uniformit de la formule \ref{superformule}. 
\par
\noindent Pour le thorme \ref{Theo1bis}, le lecteur se
reportera \`a \cite{RS08} ou \`a \cite{RS10}.\\
Dans le calcul prcdent nous avons dtermin le comportement asymptotique du
polyn\^ome pr\'edicteur
$\Phi^*_N(\chi)=\displaystyle{\sum_{k=0}^{N} \frac{(T_N(f))^{-1}_{k+1,1}}
{\sqrt{T_N(f))^{-1}_{1,1}}}}$,  on d\'eduit alors  celui de $\Phi_{N}$ via
l'op\'eration
$$
\Phi_N(z)=z^N \overline{ \Phi^*_N} ({1\over z}).
$$
\subsection{Preuve du th\'eor\`eme \ref{Theo2}}
Au lieu de considrer $\Phi_N (e^{i \theta} )$ nous considrons $\Phi_N^*(e^{i \theta})$ le polynme
prdicteur associ. La relation \ref{predi} montre que $\Phi_N (1) = \Phi_N ^*(1)$.
Posons 
$\Phi_N^*(e^{i \theta}) = \displaystyle{\sum_{l=0} ^N \beta_{N,l} (e^{i\theta})^l}.$
Ecrivons 
$\Phi_N ^*(1)= S_1+ S_2+S_3$ avec 
$ S_1= \displaystyle{\sum_{l=0} ^{[N\epsilon]} \beta_{N,l}}$,
$ S_2= \displaystyle{\sum_{l=[N\epsilon]+1} ^{N - [N\epsilon]} \beta_{N,l}}$,
$ S_3= \displaystyle{\sum_{l=N - [N\epsilon]+1 } ^{N} \beta_{N,l}}$,
avec $\epsilon$ un r\'eel positif tendant vers z\'ero..
 En utilisant le thorme \ref{Theo1} nous obtenons 
 \begin{eqnarray*}
 S_1 &=& 
   \left( \sum_{l=0} ^{[N \epsilon]}  \beta_l ^{(\alpha)} -\frac{\alpha^2} {N}  \beta_l ^{(\alpha+1)}
+ O(\frac{l}{N^2}) \right) \\
&=&  \left( \beta_{[N \epsilon]} ^{(\alpha+1)} - 
\frac{\alpha^2} {N}  \beta_{[N \epsilon]} ^{(\alpha+2)}
-\frac{[N \epsilon]}{N} O(\frac{l}{N}) \right).
\end{eqnarray*}
Comme on a d'autre part $$\beta^{(\alpha+1)} _l = \frac{1}{\Gamma (\alpha+1) c_1(1)}  l^\alpha+ O(l^\alpha) 
\quad \mathrm{et} $$
 $$\beta^{(\alpha+2)} _l = \frac{1}{\Gamma (\alpha+2)c_1(1)} l^{\alpha+1} + 
O(l^\alpha+1)$$ 
nous pouvons conclure 
$$ S_1 = \frac{1}{\Gamma (\alpha+1)c_1(1)}  [N \epsilon]^\alpha \left( 1+ \frac{\epsilon}{\alpha+1}
+ \epsilon^2 O(1)\right)$$
ou encore 
$$ S_{1}=  \frac{1}{\Gamma (\alpha+1)c_1(1)}  [N \epsilon]^\alpha + o (N^\alpha).$$
Considrons maintenant la somme $S_2$. Cette fois ci d'aprs le thorme 
\ref{Theo1bis} on a 
$$S_2  \sim  \frac{1}{\Gamma (\alpha)c_1(1)} \displaystyle{ \sum_{l = [N\epsilon]+1} ^{N- [N \epsilon]} \frac{l}{N} ^{\alpha-1} 
(1-\frac{l}{N} )^\alpha }. $$
Nous allons maintenant distinguer pour notre dmonstration les cas $\alpha >0$ et $\alpha <0$.\\
\textbf{Le cas $\frac{-1}{2} <\alpha <0$.}\\
Nous pouvons crire 
\begin{eqnarray*}
 S_2 &=& \frac{1}{\Gamma (\alpha) c_1(1)}  N^\alpha 
\left( \int_{( [N \epsilon]+1)/N } ^1 t^{\alpha-1} \left( (1-t)^\alpha -1 \right) dt \right.\\
 &+&
\left. \int_{( [N \epsilon]+1)/N } ^1 t^{\alpha-1} dt \right) + o(N^\alpha)\\
 &=& \frac{1}{\Gamma (\alpha) c_1(1)}  N^\alpha 
  \left( \int_0^1 t^{\alpha-1} \left( (1-t)^\alpha -1 \right) dt + \frac{1}{\alpha} \right) 
\\  &-&  \frac{[N \epsilon]^\alpha}{ \alpha \Gamma (\alpha) c_1(1)} +o(N^\alpha)
\end{eqnarray*}
Compte tenu de $\beta_{k} = O\left({\beta_{N-k}^{(\alpha+1)}\over N}\right) $
pour $k \in [N -  [N \epsilon]]$ on obtient \'evidemment comme plus haut 
$S_{3}= O (N^\alpha \epsilon ^{\alpha+1}) = o (N^\alpha)$. 
D'o\`u 
$$ \Phi_N ^*(1) = S_{1}+ S_{2}+ S_{3}=
 \overline{\beta_0 ^{(\alpha)} } \frac{1}{\Gamma (\alpha) c_1(1)}  N^\alpha 
  \left( \int_0^1 t^{\alpha-1} \left( (1-t)^\alpha -1 \right) dt + \frac{1}{\alpha} \right) 
+ o(N^\alpha). $$
\textbf{Le cas $0<\alpha <\frac{1}{2} $.}\\
Puisque $\alpha>0$ on a tout de suite 
$ S_{1}= O(N^\alpha)$ et d'autre part 
$$S_{2} =  \frac{1}{\Gamma (\alpha) c_1(1)}  N^\alpha 
  \left( \int_0^1 t^{\alpha-1} (1-t)^\alpha  dt  \right) + o(N^\alpha),$$
Puisque l'on a toujours $S_{3}= o(N^\alpha)$ on peut \'ecrire 
$$ \Phi_{N}^*(1) =    \frac{1}{\Gamma (\alpha) c_1(1)}  N^\alpha 
  \left( \int_0^1 t^{\alpha-1} (1-t)^\alpha  dt  \right) + o(N^\alpha).$$
  \textbf{Calcul des drives}\\
  D\'emontrons maintenant la formule dans le cas o\`u $p$ est un entier strictement positif. Avec les 
  m\^emes notations que prcdemment nous avons 
$$(\Phi_N^*)^{(p)} (e^{i \theta}) = 
\displaystyle{  \sum_{l=p} ^N l(l-1)\cdots (l-p+1)  \beta_{N,l} (e^{i\theta})^l}.$$
Nous sommes ramen\'es \`a calculer 
$(\Phi_N ^*)^{(p)}(1)= S_1+ S_2+S_3$ avec 
$ S_1= \displaystyle{\sum_{l=0} ^{[N\epsilon]} l^p \beta_{N,l}}$,
$ S_2= \displaystyle{\sum_{l=[N\epsilon]+1} ^{N - [N\epsilon]} l^p\beta_{N,l}}$,
$ S_3= \displaystyle{\sum_{l=N - [N\epsilon]+1 } ^{N} l^p \beta_{N,l}}$,
avec $\epsilon$ un r\'eel positif tendant vers z\'ero..
 En utilisant le thorme (\ref{Theo1}) nous obtenons 
 \begin{eqnarray*}
 S_1 &=& 
  \left( \sum_{l=0} ^{[N \epsilon]} l^p( \beta_l ^{(\alpha)} -\frac{\alpha^2} {N}  \beta_l ^{(\alpha+1)})
+ O(\frac{l}{N^2}) \right) \\
&=& O(N^{\alpha +p} \epsilon ^{\alpha+p})= o(N^{\alpha +p}).
\end{eqnarray*}
  On obtient de m\^eme
  $S_{3}= o(N^{\alpha +p})$, et en utilisant le th\'eor\`eme (\ref{Theo1bis})
  il est facile de v\'erifier que 
  $$ S_{2}=
 \frac{1}{\Gamma (\alpha) c_1(1)}  N^{\alpha+p}
  \left( \int_0^1 t^{\alpha+p-1} (1-t)^\alpha  dt  \right) + o(N^{\alpha+p}),$$
  ce qui ach\`eve la d\'emonstration de ce cas.\\
   \subsection{D\'emonstration du th\'eor\`eme \ref{Theo22}}
  Remarquons tout d'abord que bien sr $\Phi_{N} (1) = \overline {\Phi^*_{N}(1)}.$
  Pour le calcul des d\'eriv\'ees d'ordre sup\'erieur ou \'egal \`a un 
  nous allons l aussi devoir distinguer le cas $\alpha$ positif du cas 
  $\alpha$ ngatif. 
   L aussi trois cas sont  distinguer.\\
  Supposons d'abord $\alpha $ positif.
  Consid\'erons les sommes 
 $ \displaystyle{\sum_{u=p}^N \overline{\beta_{N-u}} u (u-1)(u-2) \cdots (u-p+1)}$
 dont la partie principale est en fait 
 \mbox{$ \displaystyle{\sum_{u=p}^N \overline{\beta_{N-u }}u^p}.$}
 L\`a aussi nous pouvons utiliser une d\'ecomposition
 $$ \sum_{u=0}^N \beta_{N-u }u^p =
 \sum_{u=0}^{[N\epsilon]} \overline{\beta_{N-u }} u^p 
 + \sum_{u=[N\epsilon]}^{N-[N\epsilon]} \overline{ \beta_{N-u }} u^p
 +\sum_{N-[N\epsilon]} ^N  \overline{\beta_{N-u }} u^p,
 $$
 o\`u $\epsilon$ est un r\'eel strictement positif qui tend vers z\'ero.
En utilisant le th\'eor\`eme \ref{Theo1} nous pouvons \'ecrire 
$$  \sum_{u=0}^{[N\epsilon]} \overline{\beta_{N-u }} u^p 
\sim \frac{\bar c_{1}(1)}{ c_{1}(1)}  \sum_{u=0}^{[N\epsilon]}
\overline{ \beta^{(\alpha+1)}}_{u} u^p, $$
 et l'on v\'erifie facilement que ce dernier terme est d'ordre 
 $O(N^{\alpha+p} \epsilon ^{\alpha+1+p}) =o(N^{\alpha+p})$.
 Nous avons de m\^eme 
 $$\sum_{N-[N\epsilon]} ^N  \overline{\beta_{N-u}}u^p \sim
\sum_{u=0}^{[N\epsilon]} u^p 
 \left( \overline{\beta_{u}^{(\alpha)}} - \frac{\alpha^2}{N}  
\overline{ \beta_{u}^{(\alpha+1)}}\right).$$
 Il est l\`a encore clair que ce terme est d'ordre $O(N^{\alpha+p}) \epsilon ^{\alpha+p})
 = o(N^{\alpha+p}).$
 Reste \`a calculer 
 $\displaystyle{ \sum_{u=[N\epsilon]}^{N-[N\epsilon]} \overline{ \beta_{N-u }u^p}}$
 ce qui peut aussi s'\'ecrire 
 $\displaystyle{ \sum_{u=[N\epsilon]}^{N-[N\epsilon]} \bar  \beta_{u }(N-u)^p}$,
  ou encore avec le th\'eor\`eme \ref{Theo1}
  $$ \sum_{u=[N\epsilon]}^{N-[N\epsilon]} \bar \beta_{u } (N-u)^p
  = \frac{1}{\Gamma (\alpha) \bar c_{1}(1)}
  \sum_{u=[N\epsilon]}^{N-[N\epsilon]} N^p 
 \left( u^{\alpha-1} (1- \frac{v}{N}) ^{\alpha+p}  +o(N^{\alpha-1})\right).$$
 La formule d'Euler et Mac-Laurin et l'hypoth\`ese $\epsilon$ tend vers z\'ero 
 permettent de conclure
 $$  \sum_{u=[N\epsilon]}^{N-[N\epsilon]} \bar  \beta_{u } (N-u)^p= 
\frac{ N^{\alpha+p} } {\Gamma (\alpha) \bar c_{1}(1)} \int _{0}^1 (1-t)^{\alpha+p} t^{\alpha-1}
dt +o(N^{\alpha+p}),$$
 et finalement 
 $$ \Phi_{N}^{(p)}(1)= \frac{ N^{\alpha+p} } {\Gamma (\alpha) \bar c_{1}(1)} \int _{0}^1 (1-t)^{\alpha+p} t^{\alpha-1}
dt +o(N^{\alpha+p}).$$
 Ce qui est l'\'enonc\'e attendu.\\
 Nous allons maintenant traiter le cas $\alpha$ n\'egatif.\\
 En introduisant la m\^eme dcomposition que dans le cas prcdent, considrons tout d'abord la somme
 $S_1= \displaystyle{ \sum_{N - N \epsilon}^N \overline{\beta_{N-u} }u^p }$.
 On a 
 \begin {eqnarray*}
 S_1&=& \sum_{u= N - [N \epsilon]}^N 
 \left(\overline{ \beta^{(\alpha)} _{N-u}}  - \frac{\alpha^2} {N}  \overline{ \beta^{(\alpha+1)} _{N-u}}\right) u^p
 =
  \sum_{u= 0}^{[N \epsilon]} 
 \left( \overline{\beta^{(\alpha)} _{u}}  - \frac{\alpha^2} {N}  \overline{ \beta^{(\alpha+1)} _u }\right) (N-u)^p\\
 &=& \sum_{k=0}^p (-1) ^k  C_p^k N^{p-k} 
 \left(  \sum_{u= 0}^{[N \epsilon]}  
  \left( \overline{\beta^{(\alpha)} _{u}}  - \frac{\alpha^2} {N}  \overline{ \beta^{(\alpha+1)} _{u} }\right) u^k\right)\\
  &\sim & 
 \sum_{k=0}^p  (-1) ^k  C_p^k N^{p-k}  
 \left( \frac{(N \epsilon )^{\alpha+k} } {(\alpha +k) \Gamma (\alpha) } 
 - \frac{\alpha^2 (N \epsilon )^{\alpha+k+1} } {(\alpha +k+1) \Gamma (\alpha+1) } \right) \frac{1}{\bar c_{1}(1)}\\
 &=& \frac{N^{p+\alpha} \epsilon ^\alpha}{\Gamma (\alpha+1)\bar c_{1}(1)}+ \epsilon O(N^{p+\alpha}\epsilon ^\alpha)
 = \frac{N^{p+\alpha} \epsilon ^\alpha}{\Gamma (\alpha+1)\bar c_{1}(1)} + o(N^{p+\alpha}\epsilon ^\alpha).
 \end{eqnarray*}
 On a d'autre part , si l'on pose 
 $\displaystyle{S_2 = \sum_{u=[N \epsilon]}^{N - [N \epsilon]}\overline{ \beta_{N-u} } u^p}$
 \begin{eqnarray*}
 S_2 &=& \sum_{u=[N \epsilon]}^{N - [N \epsilon]}\overline{ \beta_{N-u}} u^p
 =  \sum_{u=[N \epsilon]}^{N - [N \epsilon]} \bar \beta_{u} (N-u)^p\\
 &=& \frac {N^{p+\alpha-1}}{\Gamma (\alpha)\bar c_{1}(1)} \sum_{u=[N \epsilon]}^{N - [N \epsilon]}  
  \left( u ^{\alpha-1} (1-\frac{u}{N}) ^{\alpha+p}\right) \left(1+o(1)\right)\\
  &=& \frac {N^{p+\alpha-1}}{\Gamma (\alpha)\bar c_{1}(1)} 
\int_\epsilon ^{1-\epsilon} t^{\alpha-1} (1-t)^{\alpha+p}dt \\
&=&  \frac {N^{p+\alpha-1}}{\Gamma (\alpha)\bar c_{1}(1)} 
\int_\epsilon ^{1-\epsilon} t^{\alpha-1} ((1-t)^{\alpha+p}-1)dt 
 + \frac {N^{p+\alpha-1}}{\Gamma (\alpha)\bar c_{1}(1)} \int_\epsilon ^{1-\epsilon} t^{\alpha-1} dt 
 \\ &=&  \frac {N^{p+\alpha-1}}{\Gamma (\alpha)\bar c_{1}(1)} 
\int_\epsilon ^{1-\epsilon} t^{\alpha-1} ((1-t)^{\alpha+p}-1)dt 
 -  \frac{N^{p+\alpha} \epsilon ^\alpha}{\Gamma (\alpha+1)\bar c_{1}(1)} + \frac{N^{p+\alpha}}{\Gamma (\alpha+1)\bar c_{1}(1)} 
 \end{eqnarray*}
 La troisime somme intervenant dans la dcomposition tant clairement ngligeable, on a
 finalement 
 $$ \Phi_N^{(p)} (1) =  \frac {N^{p+\alpha-1}}{\Gamma (\alpha)\bar c_{1}(1)} 
\left( 
\int_\epsilon ^{1-\epsilon} t^{\alpha-1} ((1-t)^{\alpha+p}-1)dt + \frac{1}{\alpha} \right)
+ o (N^{p+\alpha-1}),
$$
 ce qui donne imm\'ediatement le r\'esultat. 
 
    \subsection{Preuve du th\'eor\`eme \ref{proba1}.}
  \subsubsection{D\'eterminant des op\'erateurs de Fredholm dans $L^2(\mathbb T)$}
  Pour $(\theta, \theta' )$ dans $]-\pi, \pi] \times  ]-\pi, \pi]$ nous posons
  $$ k(\theta, \theta') =
  \sum_{m=1}^N \phi_{m}(\theta) f_{m} (\theta')$$
  o\`u pour tout entier $m$ tles fonctions $\phi_{m}$ et $f_{m}$ sont dans
  $L^2(\mathbb T)$. 
 Si $X$ est un intervalle $]-\pi, \pi]$ on peut alors d\'efinir 
  un op\'erateur $K$ sur $L^2 (X)$ en posant, pour tout $x$ dans 
  $L^2 (X)$, 
  $$ K_{x} (\theta) = \int_{X} k (\theta,\theta') x(\theta') d \theta'.$$
  Il est alors connu (voir, par exemple,  \cite {Gohgolkru} ou \cite{ RS7580} )
   que $\det (I+K)$ est la quantit\'e d\'efinie par 
  $$
  \det(I+K) = 1+ \sum_{m=1}^N \frac{1}{m!} \int _{X^m} 
  \det (k(\theta_{i}, \theta_{j})_{i,j=1,\cdots,m} d\theta_{1}\cdots d \theta_m
  $$
  ou encore en posant que 
 l'intgrale sur 
 $X^{0}$ vaut 1 on peut crire, de manire plus synthtique
  $$
  \det(I+K) =  \sum_{m=0}^N \frac{1}{m!} \int _{X^m} 
  \det (k(\theta_{i}, \theta_{j})_{i,j=1,\cdots,m} d\theta_{1}\cdots d \theta_m
  $$

   De m\^eme on peut d\'efinir la trace de $K$ par 
  $$ \Tr K = \int_{X} k(\theta, \theta) d\theta.$$
  D'autre part si $(K_{n})_{n \in \mathbb N}$ est une suite d'op\'erateurs ainsi 
  d\'efinis qui converge vers un op\'erateur $A$  au sens des op\'erateurs sur 
  $L^2(X)$, on peut poser 
  $$ \det (Id+A) = \lim_{n \rightarrow \infty} \det (I+K_{n}) 
  \quad \mathrm {et} \quad  \Tr A = \lim_{n \rightarrow \infty} \Tr (K_{n}).
  $$
  On sait que $A$ est alors un op\'erateur fix\'e (voir \cite{Gohgolkru})
  la fonction $\lambda \rightarrow \det (I+\lambda A)$ est une fonction enti\`ere sur 
  $\mathbb C$ 
  dont le d\'eveloppement en s\'erie enti\`ere est donn\'e par 
  $$ \det (Id +\lambda A) = 1 + \sum_{n=1}^\infty \frac{c_{n}(A)} { n!}  \lambda^n.$$
  Les coefficients $c_{n}$ \'etant d\'efinis par 
  
  $$
  c_{n}(A) = \det \left( 
  \begin {array}{cccccc}
  \Tr (A) & n-1 & 0 &\cdots & 0 & 0\\
  \Tr (A^2) & \Tr (A) & n-2 & \cdots &0 &0 \\
  . & .&.&.&.&.\\
   . & .&.&.&.&.\\
    . & .&.&.&.&.\\
    \Tr (A^{n-1}) & \Tr (A^{n-2}) &  \Tr ( A^{n-3}) & \cdots & \Tr (A) & 1\\
     \Tr (A^{n}) & \Tr (A^{n-1}) &  \Tr ( A^{n-2}) & \cdots & \Tr (A^2) & \Tr (A)\\
  \end{array}
  \right)
  $$
 Ce qui peut s'\'ecrire aussi 
 $$ 
 \det (Id + \lambda A) = 
 \exp \left ( \sum _{m=1}^\infty \frac{(-1) ^{m+1}} { m!} \Tr (A^m)\lambda^m\right).
 $$
 On peut alors \'enoncer le r\'esultat :
\begin{prop}
 L'op\'erateur $Id +A$ est inversible si et seulement si  $ \det (Id + A) \not=0$.  
  \end{prop}
  \subsubsection {Polyn\^omes orthogonaux et fonction de corr\'elation}
  Consid\'erons l'\'egalit\'e 
  $$ \Phi^*_{N} (e^{i \theta}) 
= \sum _{j=0}^N \frac{(e^{i \theta}-1)^j }{j!}  (\Phi^*_{N})^{(j)} (1). $$
Nous pouvons \'ecrire, pour un entier $k_{0}$ fix\'e 
$$ \Phi^*_{N} (e^{i \theta}) 
= \sum _{j=0}^{k_{0}} \frac{(e^{i \theta}-1)^j }{j!}  (\Phi^*_{N})^{(j)} (1)
 +  N^{\alpha} O\left( \frac{ \vert \theta ^{k_0+1}\vert } { (k_{0}+1)!}\right ). $$
 En effet la dmonstration du th\'eor\`eme \ref{Theo2} permet d'obtenir , si $j>k_{0}$, et en 
 remarquant que $ l (l-1)(l-2) \cdots (l-p+1) <l^p$ :
$$ \vert  \frac{(e^{i \theta/N}-1)^j }{j!}  (\Phi^*_{N})^{(j)} (1) \vert \le  \frac{N ^\alpha \vert  i \theta)^j \vert}
{c_1(1) \Gamma(\alpha)} \left ( \int_{0}^1 \frac{x^{\alpha +j-1} (1-x)^\alpha}{j!} dx +o(N^\alpha) \right)$$
uniform\'ement par rapport \`a $j$.
Ce qui donne, toujours en utilisant le thorme \ref{Theo2}
$$
\Gamma (\alpha) c(1)  \Phi^*_{N} (e^{i \theta/N})  = N^\alpha\left(
\sum_{j=0}^{k_{0}} \int_{0}^1 \frac {x^{\alpha+j-1}(i \theta)^j }{j!} (1-x)^\alpha dx +O\left( \frac{ \vert \theta ^{k_0+1}\vert } { (k_{0}+1)!}\right )
+o(1)\right).
$$
Ceci peut encore s'crire 
\begin{equation} \label{un}
\Gamma (\alpha) c_{1}(1)   \Phi^*_{N} (e^{i \theta/N})  = {N^\alpha}
\left(\int_{0}^1 e^{i\theta x} x^{\alpha-1} (1-x)^\alpha dx + R(k_{0})\right).
\end{equation}
avec 
$$ R(k_{0}) =- \sum_{k_0+1} ^{+ \infty}  \int_{0}^1 \frac {x^{\alpha+j-1}(i \theta)^j }{j!} (1-x)^\alpha dx 
+O\left( \frac{ \vert \theta ^{k_0+1}\vert } { (k_{0}+1)!}\right ).
$$
Si maintenant nous nous fixons un r\'eel $\epsilon$ positif, il est clair que l'on peut choisir le r\'eel $k_{0}$ tel que 
$\vert R(k_{0}) \vert <\epsilon $. Nous pouvons donc finalement \'ecrire 
\begin{equation} \label{unono}
\Gamma (\alpha) c_{1}(1)   \Phi^*_{N} (e^{i \theta/N})  = {N^\alpha}
\left(\int_{0}^1 e^{i\theta x} x^{\alpha-1} (1-x)^\alpha dx +o(1)\right).
\end{equation}
On obtient de m\^eme que 
\begin{equation}\label {deux}
\Gamma (\alpha) c_{1}(1)  \Phi_{N} (e^{i \theta/N})  = {N^\alpha}
\left( \int_{0}^1 e^{i \theta(1-x)} x^{\alpha-1} (1-x)^{\alpha} dx +o(1)\right).
\end{equation}
Introduisons maintenant le noyau de Christofel-Darboux d'ordre $N$, \`a savoir la 
quantit\'e

$$ K_N (e^{i\theta}, e^{i\theta'}) =  \sqrt{f(e^{i \theta}) }
\sqrt {f(e^{i \theta'})}
\sum_{m=0}^{N-1} \frac{1}{h_m} \overline{ \Phi_m (e^{i \theta}) }
\Phi_m (e^{i \theta'}) 
$$ 
o\`u les $h_{m}$ sont des constantes de normalisations d\'efinies par 
$$
\int_{- \pi} ^{\pi} f(e^{i \theta}) \overline{ \Phi_m (e^{i \theta}) } \Phi_n (e^{i \theta}) d\theta 
= h_m \delta_{m,n}.
$$ 
On sait que l'on peut \'ecrire alors, si $\theta\not= \theta'$  
(voir \cite{SZEG})
$$  K_N (e^{i\theta}, e^{i\theta'})=  \frac{1}{h_{N}} \sqrt{f(e^{i \theta}) }
\sqrt {f(e^{i \theta'})}
\frac {\overline{\Phi_{N} ^* (e^{i\theta})} \Phi_{N} ^* (e^{i\theta'})
-\overline {\Phi_{N}  (e^{i\theta}) } \Phi_{N} (e^{i\theta'})}
{ (1-e^{i(\theta'-\theta)} )}. $$

En posant dans la formule ci-dessus $\theta = \frac{u}{N}, \theta = \frac{v}{N}$
on obtient, en utilisant les formules (\ref{un}) et (\ref{deux}) 
\begin{equation} \label{nepasnegliger}
  \Gamma^2 (\alpha)  c_{1} (1) K_{N} (e^{i u/N},e^{iv/N} ) 
=\frac{1}{h_N}  N^{2 \alpha} \sqrt {f(e^{i u/N}) }
\sqrt{ f(e^{iv/N})} \frac{K(u,v)+o(N^{2\alpha})}
{ 1-e^{i( u/N- v/N )}}
\end{equation}
avec si $u\not=v$
\begin{eqnarray*}
K(u,v)&= &
  \int_{0}^1 e^{-iux} x^{\alpha-1} (1-x)^\alpha dx
 \int_{0}^1 e^{ivx} x^{\alpha-1} (1-x)^\alpha dx -\\
& -& \int_{0}^1 e^{-iu(1-x)} x^{\alpha-1} (1- x)^{\alpha} dx 
\int_{0}^1 e^{iv(1-x)} x^{\alpha-1} (1-x)^{\alpha} dx 
\end{eqnarray*}
 Ou en  posant 
$$
\psi (\alpha,u) = \int_0 ^{1} x^{\alpha-1} (1-x) ^\alpha e^{i u x} dx,
$$
$$ 
K(u,v)=\psi (\alpha,-u) \psi(\alpha,v) 
- \psi (\alpha,u) \psi(\alpha,-v) e^{i(v-u)}
$$ 
et en remplaant $\sqrt f$ et $(1-e^{i(u/N-v/N)} )$ par un quivalent, il vient,
toujours si $u\not=v$, 
 $$ \Gamma ^2 (\alpha)  c_{1}(1)  K_{N} (e^{i u/N},e^{iv/N} ) 
= \frac{1}{h_N} N (uv)^\alpha \frac  {K(u,v)} {i(u-v)}+ o(N).$$
Pour conna\^{i}tre la valeur $u=v$ consid\'erons la d\'eriv\'ee 
au point $v=u$ de la fonction $G_{u}$ d\'efinie par 
$$ G_{u}(v) = e^{iu} \psi (\alpha,-u) \psi (\alpha,v) 
- e^{iv}  \psi (\alpha,u) \psi (\alpha,-v).$$
La limite quand $v \rightarrow u$ de $K_{N} (e^{i u/N},e^{iv/N} )$
est alors $ \frac{1}{h_N} N (u)^{2\alpha}\left( \frac{e^{-iu} }{-i} \right)G'_u(u).$
 Ce qui donne, en posant 
 $\tau(\alpha,u) = \int_{0}^1 x^\alpha (1-x)^\alpha e^{iux} dx$,
 $$ \Gamma ^2 (\alpha) c_1 (1) K_N (e^{iu/N}, e^{iu/N}) = \frac{1}{h_N} 
 N u^{2\alpha} \left( -2 \Re \left(\psi (\alpha,u) \tau(\alpha,-u)\right  ) 
 +\vert \psi (\alpha,u) \vert ^2 \right) +o(N).$$
On sait d'autre part (voir \cite{Ld}) qu'avec les notations du corollaire \ref{uno} 
$$h_{m}= \int_{-\pi}^{\pi} f(e^{i \theta}) \vert \Phi_{N}(e^{i \theta})\vert^2 d \theta,$$
c'est \`a dire 
$$h_{m}= \left(T_{m}(f)\right)^{-1}_{1,1}. $$ Nous avons d'autre part obtenu dans 
un pr\'ec\'edent travail ( \cite{RS08},\cite{RS10}) que 
$$ \vert \beta_{0}^\alpha\vert ^2(1-Ê\frac{\alpha^2}{N}) \sim  \left(T_{m}(f)\right)^{-1}_{1,1} .$$
Finalement l'hypoth\`ese   $ \beta_{0}^\alpha=1$ nous donne,   
avec $u\not=v$,
$$ \Gamma^2 (\alpha) c_{1} (1)  K_{N} (e^{i u/N},e^{iv/N} ) 
=  N (uv)^\alpha \frac  {K(u,v)} {i(u-v)}+ o(N),$$
et si $u=v$ 
 $$ \Gamma ^2 (\alpha) c_1 (1) K_N (e^{iu/N}, e^{iu/N}) =
 N u^{2\alpha} \left( 2 \Re \left(\psi (\alpha,u) \tau(\alpha,u)\right  ) 
 +\vert \psi (\alpha,u) \vert ^2 \right) +o(N).$$

Si maintenant on d\'esigne par $R_{k}^{(N)}$ la fonction de corr\'elation d\'efinie par 
$$ R_{k}^{(N)} (\theta_{1}, \theta_{2}, \cdots, \theta_{k}) = \frac{1}{(N-k)!} 
\int_{]-\pi, \pi[ ^{N-k}} P_{N}(\theta_{1},\cdots, \theta_{N}) d\theta_{k+1}\cdots 
d\theta_{N}$$
avec 
$$ P_{N}(\theta_{1},\cdots, \theta_{N}) = \prod_{j=1}^N 
\vert 1 - e^{i \theta_{j}} \vert^{2 \alpha} c_{1} (\theta_{j})
\prod _{j<i} \vert e^{i \theta_{j}} - e^{i \theta_{i}} \vert ^2.$$

On sait qu'on a alors la relation fondamentale 
\begin{prop}  \label{fondamentale}
pour tout intervalle $I$ contenu dans $] - \pi , \pi[$ et tout entier 
$m \in \{0, \cdots,N\}$ on a 
\begin{equation} 
P \left( \Big \vert \{ i \le N, \quad \theta_{i} \in I \} \Big \vert  =m \right) =
\frac{(-1)^m }{m} \sum_{k=m}^N \frac{(-1)^k} {(k-m) !} \int _{I^k}
R_{k}^{(N)} (\theta) d^k (\theta).
\end{equation}
\end{prop}
Des calculs classiques sur les polyn\^omes orthogonaux donnent la relation
$$ R_{k}^{(N)} (\theta_{1} \cdots \theta_{k})
= \det [ K_{N}(\theta_{i},\theta_{j}) _{i,j =1 \cdots k}] .$$
C'est \`a dire que 
\begin{equation} \label{fondamentale2}
P \left( \Big \vert \{ i \le N, \quad \theta_{i} \in I \} \Big \vert =m\ \right)=
\frac{(-1)^m }{m} \sum_{k=m}^N \frac{(-1)^k} {(k-m) !} \int _{I^k}
 \det [ K_{N}(\theta_{i},\theta_{j}] _{i,j =1 \cdots k}] d^k (\theta).
\end{equation}
On obtient donc finalement, en posant maintenant $I =[u,v]$,
$-\pi<u<v<\pi$ et $I_{N} =[\frac{u}{N}, \frac{v}{N}],$ la formule 
\begin{equation} \label{fondamentale3}
P \left( \Big \vert \{ i \le N, \quad \theta_{i} \in I_{N} \} \Big \vert  =m \right) =
\frac{(-1)^m }{m} \sum_{k=m}^N \frac{(-1)^k} {(k-m) !} \int _{I_{N}^k}
\det [ K_{N}(\theta_{i},\theta_{j}] _{i,j =1 \cdots k}]  d^k (\theta).
\end{equation}
En posant le changement de variables $u_{i} = N \theta_{i}$ on obtient
\begin{align} \label{fondamentale4}
& P \left( \Big \vert \{ i \le N, \quad   \theta_{i} \in I_{N} \} \Big \vert = m   \right) =\\
&= \frac{(-1)^m }{m} \sum_{k=m}^N \frac{(-1)^k} {(k-m) !} \int _{I^k}
N^{-k}  \det [ K_{N}(u_{i}/N,u_{j}/N] _{i,j =1 \cdots k}]  d^k (u), \nonumber
\end{align}
ou encore 
\begin{equation} \label{fondamentale5}
P \left( \Big \vert \{ i \le N, \quad \theta_{i} \in I_{N} \} \Big \vert  =m \right) =
\frac{{-1}^m }{m} \sum_{k=m}^N \frac{(-1)^k} {(k-m) !} \int _{I^k}
  \det [ H_{N}(u_{i},u_{j}] _{i,j =1 \cdots k}]   d^k (u),
\end{equation}
avec,si $u\not=v$ 
$$ H_N (u,v) = \frac{1} {\Gamma^2 (\alpha) c_1(1)}
 (uv)^\alpha \frac  {K(u,v)} {i(u-v)}+ o(1),$$
et si $u=v$ 
$$ H_N (u,u) = \frac{1} {\Gamma^2 (\alpha) c_1 (1)}
 u^{2\alpha} \left( -2 \Re \left(\psi (\alpha,u) \tau(\alpha,-u)\right  ) 
 +\vert \psi (\alpha,u) \vert ^2 \right)+o(1).$$
 En revenant \`a la d\'efinition des d\'eterminants des op\'erateurs de Fredholm
nous pouvons \'ecrire
 \begin{equation} \label{fonda2}
\sum_{k=0}^N \frac{(-1)^k} {k !} \int _{I^k}
 \det [ H_{N}(u_{i},u_{j}) _{i,j =1 \cdots k}]   d^k (u),
= \det [(I-  H_{N}Ê )\vert_{L^2(I)}]
\end{equation}
Ce qui donne, au moyen d'un changement de variables, 
 \begin{equation} \label{fonda}
 \sum_{k=0}^N \frac{(-\gamma )^k} {k !} \int _{I^k}
 \det [ H_{N}(u_{i},u_{j}) _{i,j =1 \cdots k}]  d^k (u),
= \det[ (I-\gamma   H_{N} )\vert_{L^2(I)}]
\end{equation}
ou encore 
\begin{equation}
 \sum_{k=m}^N \frac{(-1)^k \gamma^{k-m} } {(k-m) !} \int _{I^k}
 \Big \vert \det [ H_{N}(u_{i},u_{j}) _{i,j =1 \cdots k}] \Big \vert d^k (u)=
\frac{d^m}{d\gamma^m} \det [(I-\gamma  H_{N})\vert_{L^2(I)}].
 \end{equation}
soit 
\begin{equation}
\sum_{k=m}^N \frac{(-1)^k } {(k-m) !} \int _{I^k}
\det [ \Big \vert H_{N}(u_{i},u_{j}) _{i,j =1 \cdots k}] \Big \vert d^k (u)=
\left(\frac{d^m}{d\gamma^m} \det [(I-\gamma  H_{N} )\vert_{L^2(I)}]
\right)
\Big \vert _{\gamma=1}.
\end{equation}
Et par passage  la limite on obtient finalement
\begin{equation} \label{fundamentale6}
\lim_{N \rightarrow + \infty}P \left( \Big \vert \{ i \le N \theta_{i} \in I_{N} \} \Big \vert \right)=m\} =
 \frac{(-1)^m} {m!} \left( \frac{d}{d\gamma}\right)^m 
\det [ (Id - \gamma \mathcal K)\vert _{L^2(I)}] \Big \vert _{\gamma =1}.
\end{equation}
Ceci termine la d\'emonstration du th\'eor\`eme \ref{proba1}. Le
deuxime point du thorme est alors la traduction du cas $m=0$.\\
La d\'emonstration du th\'eor\`eme \ref{proba2} repose bien s\^ur sur les m\^emes
id\'ees que celle ci.

\subsection{D\'emonstration du th\'eor\`eme \ref{proba2}}
Nous pouvons de nouveau \'ecrire
$$ \Phi^*_{N} (e^{i \theta}) 
= \sum _{j=0}^N \frac{(e^{i \theta}-1)^j }{j!}  (\Phi^*_{N})^{(j)} (1).$$
Comme dans le cas $\alpha$ positif le th\'eor\`eme (\ref{Theo2}) permet donne alors pour $j\ge 1$ :
$$  \frac{(e^{i \theta/N}-1)^j }{j!}  (\Phi^*_{N})^{(j)} (1)= \frac{N ^\alpha( i \theta)^j}
{c(1) \Gamma(\alpha)} \int_{0}^1\frac{ x^{\alpha +j-1} (1-x)^\alpha}{j!} dx +o(N^\alpha).$$
Ce qui donne finalement, en utilisant les m\^emes arguments que dans la d\'emonstration pr\'ec\'edente, 
\begin{eqnarray*}
\frac{\Gamma (\alpha) c_{1}(1) } {N^\alpha}  \Phi^*_{N} (e^{i \theta/N}) &= &
\sum_{j=1}^N \int_{0}^1 \frac {x^{\alpha+j-1}(i \theta)^j }{j!} (1-x)^\alpha dx \\
&+& \int _0^1x^{\alpha -1} \left((1-x)^\alpha -1 \right) dx + \frac{1}{\alpha} 
+o(1).
\end{eqnarray*}
Tout ceci permet de conclure que 
\begin{align}  \label{unprime}
\Gamma (\alpha) c_{1}(1)   \Phi^*_{N} (e^{i \theta/N})  &= {N^\alpha}
\left(\int_{0}^1  x^{\alpha-1} (e^{i\theta x} -1)(1-x)^\alpha dx+ \right. 
\\
 \nonumber &  +\left.
\int_{0}^1  x^{\alpha-1} ((1-x)^\alpha-1) dx+ \frac{1}{\alpha}
+o(1)\right)
\end{align}
ou encore 
\begin{equation} \label{unbisbis}
\Gamma (\alpha) c_{1}(1)   \Phi^*_{N} (e^{i \theta/N})  = {N^\alpha}
\left(\int_{0}^1  x^{\alpha-1} \left(e^{i\theta x} (1-x)^\alpha -1\right) dx +
 \frac{1}{\alpha} +o(1)\right).
\end{equation}

On obtient de m\^eme, en utilisant la relation$ \Phi^*_N (z) =z^N \bar \Phi_N (\frac{1}{z}) $ que 
\begin{equation}\label {deuxprime}
\Gamma (\alpha) c_{1}(1)  \Phi_{N} (e^{i \theta/N})  = {N^\alpha} e^{i \theta}
\left( \int_{0}^1 (e^{-i \theta x}  (1-x)^{\alpha}-1 ) x^{\alpha-1} dx + 
\frac{1}{\alpha} +o(1)\right).
\end{equation}
C'est \`a dire que 
$$\frac{ \Gamma (\alpha) c_{1}(1)  \Phi_{N} (e^{i \theta/N})}{N^\alpha}
= e^{i \theta} \frac{ \Gamma (\alpha) c_{1}(1)  \Phi_{N}^* (e^{-i \theta/N})}{N^\alpha}
+o(1).
$$
 D\'esignons  maintenant par $\tilde K_{N}$ le noyau de Christoffel-Darboux d'ordre $N$.
On pose de nouveau  $\theta=\frac{u}{N}, et \theta'=\frac{v}{N}$,
 et en on utilise les formules (\ref{unprime}) et (\ref{deuxprime}). 
On obtient alors,par les m\^emes m\'ethodes que pr\'ec\'edemment, et pour $u\not=v$
\begin{equation} \label{nepasnegliger2}
  \Gamma^2 (\alpha) \vert c_{1} (1)\vert ^2 \tilde K_{N} (e^{i u/N},e^{iv/N} ) 
=\frac{1}{h_{N}}  N^{2 \alpha} \sqrt{f(e^{i u/N}) }
\sqrt { f (e^{iv/N})} \frac{\tilde K(u,v)+o(N^{2\alpha})}
{ 1-e^{i(u/N-v/N)}}
\end{equation}
avec 
$$
\tilde K(u,v) = \tilde \psi (\alpha,-u) \tilde \psi(\alpha,v) 
- \tilde \psi (\alpha,u) \tilde \psi(\alpha,-v) 
 e^{i(v-u)},$$
et
$$
\tilde \psi (\alpha,u) = \int_0 ^{1} x^{\alpha-1}( (1-x) ^\alpha  e^{-i u x}-1) dx
+ \frac{1}{\alpha}.$$
En remplaant $\sqrt f$ et $(1-e^{i(u/N-v/N)} )$ par un quivalent, il vient 
si $u\not=v$
 $$ \Gamma ^2 (\alpha)  c_{1}(1)  \tilde K_{N}
  (e^{i u/N},e^{i v/N} ) 
= \frac{1}{h_{N}} N (- uv)^\alpha \frac  {\tilde K(u,v)} {i(u-v)}
+ o(N).$$
La valeur pour $u=v$ s'obtient par un simple passage \`a la limite.
Le reste de la d\'emonstration est alors identique \`a celle du th\'eor\`eme \ref{proba1}.
\subsection{D\'emonstration du th\'eor\`eme \ref{proba3}}
Nous pouvons de nouveau \'ecrire pour $p>1$ 
$$ \Phi^*_{N} (e^{i \theta/N^p}) 
= \sum _{j=0}^N \frac{(e^{i \theta/N^p}-1)^j }{j!}  (\Phi^*_{N})^{(j)} (1).$$
L\`a aussi le th\'eor\`eme \ref{Theo2} permet d'obtenir :
$$  \frac{(e^{i \theta/N^p}-1)^j }{j!}  (\Phi^*_{N})^{(j)} (1)=  \frac{1}{N^{j(p-1)}}
 \frac{N ^\alpha ( i \theta)^j}
{c_{1}(1) \Gamma(\alpha)} \int_{0}^1 \frac{x^{\alpha +j-1} (1-x)^\alpha}{j!} dx +o(N^\alpha).$$
Ce qui donne finalement 
$$
\frac{\Gamma (\alpha) c_{1}(1) } {N^\alpha} \Phi^*_{N} (e^{i \theta/N^p})  =
\sum_{j=0}^N \int_{0}^1 \frac {x^{\alpha+j-1} (i \theta/N^{p-1} ) ^j}{j!} (1-x)^\alpha dx 
+o(1).
$$
Ce qui peut aussi s'\'ecrire 
\begin{equation} \label{one}
\Gamma (\alpha) c_{1}(1)   \Phi^*_{N} (e^{i \theta/N^p})  = {N^\alpha}
\left(\int_{0}^1 e^{i\theta x/N^{p-1}} x^{\alpha-1} (1-x)^\alpha dx +o(1)\right).
\end{equation}
On obtient de m\^eme que 
\begin{equation}\label {two}
\Gamma (\alpha) c_{1}(1)  \Phi_{N} (e^{i \theta/N})  = {N^\alpha}
\left( \int_{0}^1 e^{i \theta(1-x)/ N^{p-1}} x^{\alpha-1} (1-x)^{\alpha} dx +o(1)\right).
\end{equation}
C'est \`a dire que 
\begin{align} \label{nepasnepasnegliger}
  \Gamma^2 (\alpha)  c_{1} (1) K_{N} (e^{i u/N^p},e^{iv/N^p} ) 
& =\frac{1}{h_N}  N^{2 \alpha} \sqrt {f(e^{i u/N^p}) } \sqrt{ f(e^{iv/N^p})}\\
& \times  \frac{K(u/N^{p-1},v/N^{p-1})+o(N^{2\alpha})}
{ 1-e^{i( u/N^p- v/N^p )}}.\nonumber
\end{align}
Ce qui permet de calculer les probabilit\'es 
$P \left( \Big \vert \{ i \le N, \quad \theta_{i} \in I_{N^p} \} \Big \vert \right)=m\}$
pour tous les entiers $m,p \ge 1$.
En posant $I =[u,v], \quad 
-\pi<u<v<\pi$ et $I_{N^p} =[\frac{u}{N^p}, \frac{v}{N^p}],$ on a\begin{equation} \label{fondamentale33}
P \left( \Big \vert \{ i \le N, \quad \theta_{i} \in I_{N^p} \} \Big \vert=m \right) =
\frac{(-1)^m }{m} \sum_{k=m}^N \frac{(-1)^k} {(k-m) !} \int _{I_{N}^k}
  \det [ K_{N}(\theta_{i},\theta_{j}] _{i,j =1 \cdots k}]  d^k (\theta).
\end{equation}
En posant le changement de variables $u_{i} = N^p \theta_{i}$ on obtient
\begin{eqnarray*} 
& &P \left( \Big \vert \{ i \le N, \quad \theta_{i} \in I_{N} \} \Big \vert  =m\right) =\\
&=& \frac{(-1)^m }{m} \sum_{k=m}^N \frac{(-1)^k} {(k-m) !} \int _{I^k}
N^{-p k}  \det [ K_{N}(u_{i}/N^p,u_{j}/N^p] _{i,j =1 \cdots k}]  d^k (u),
\end{eqnarray*}
ou encore 
\begin{eqnarray*} 
& & P \left( \Big \vert \{ i \le N, \quad \theta_{i} \in I_{N} \} \Big \vert  =m \right)=\\
& =& \frac{{-1}^m }{m} \sum_{k=m}^N  N^{k(1-p)}
\frac{(-1)^k} {(k-m) !} \int _{I^k}
\det [ H_{N}(u_{i}/N^{p-1},u_{j}/N^{p-1})_{i,j =1 \cdots k}] d^k (u),
\end{eqnarray*}
Puisque pour tout entier naturel $k\in [1,N]$
la quantit\'e $\det [ H_{N}(u_{i}/N^{p-1},u_{j}/N^{p-1})_{i,j =1 \cdots k}]$ 
tend vers z\'ero quand $N$ tend vers l'infini  par continuit de la fonction $H_N$
 il est clair que le r\'esultat annonc\'e dans l'\'enonc\'e est acquis. 
 \section{Appendice, r\'esultats num\'eriques}
         Dans chaque tableau la premi\`ere colonne donne $\Phi_{N} (1)$ la deuxi\`eme colonne donne 
 $A N^{-d}$ o\`u $A$ est une constante induite par les r\'esultats des 
th\'eor\`emes \ref{Theo2} et \ref{Theo22}. Les valeurs de $N$ \'etudi\'es sont entre $400$ 
et $640$. Rappelons que par rapport au corps de l'article $d = -\alpha$.
 \subsubsection{d=-0.2750} 
  \begin{tabular} [t] {c c}
  $\Phi_{N}(1)$ & $ A N^{-d}$\\
  11.7692 &  11.7688 \\
      11.7746 & 11.7742 \\
   11.7800 &   11.7796\\
   11.7854 & 11.7850 \\
   11.7907 &  11.7904 \\
   11.7961 & 11.7958 \\
   11.8015 &  11.8011 \\
   11.8068 &  11.8065 \\
   11.8122 &   11.8118\\
   11.8175 & 11.8172\\
  1.8228 &    11.8225\\
      \end{tabular} 
       \begin{tabular} [t] {c c}
        $\Phi_{N}(1)$ & $A N^{-d}$\\
         1   11.8282 &    11.8279\\
   11.8335 &11.8332\\
   11.8388 &  11.8385\\
   11.8441 & 11.8438\\
   11.8494 &   11.8491\\
   11.8547 &   11.8544\\
   11.8600 &11.8597\\
   11.8653 & 11.8650\\
   11.8705 &  11.8703\\
           11.8758 &  11.8756 \\
        \end{tabular} 
       \begin{tabular} [t] {c c}
        $\Phi_{N}(1)$ & $A N^{-d}$\\
   11.8811 & 11.8809 \\
   11.8863 &   11.8861\\
   11.8916&   11.8914\\
   11.8968&    11.8967\\
   11.9021&    11.9019\\
   11.9073&    11.9071\\
   11.9125&    11.9124\\
   11.9177&  11.9176\\
   11.9229&   11.9228\\
     11.9282&   11.9281\\
                \end{tabular}
       \begin{tabular} [t] {c c}
        $\Phi_{N}(1)$ & $A N^{-d}$\\       
   11.9334&   11.9333\\
   11.9386&    11.9385\\
   11.9437&    11.9437\\
   11.9489&    11.9489\\
   11.9541&    11.9541\\
   11.9593&  11.9592\\
   11.9644& 11.9644\\
   11.9696&   11.9696\\
   11.9748&  11.9747\\
   11.9799&  11.9799\\

        \end{tabular}
             
       \subsubsection {d= -0.1500}
  \begin{tabular} [t] {c c}
  $\Phi_{N}(1)$ & $A N^{-d}$\\
    5.4235 &5.4234\\
    5.4249&  5.4247\\
    5.4262&  5.4261\\
    5.4276& 5.4274\\
    5.4289&  5.4288\\
    5.4303&5.4301\\
    5.4316&    5.4315\\
    5.4329&5.4328\\
    5.4343&   5.4342\\
    5.4356&   5.4355\\
    5.4370&    5.4368\\
    5.4383& 5.4382\\
    5.4396&  5.4395\\
    5.4410&   5.4408\\
   5.4423&    5.4422\\
    5.4436& 5.4435\\
     \end{tabular} 
  \begin{tabular} [t] {c c}
  $\Phi_{N}(1)$ & $A N^{-d}$\\
    5.4449& 5.4448\\
    5.4463&  5.4462\\
    5.4476&  5.4475\\
  5.4489&   5.4475\\
5.4383& 5.4382\\
    5.4396&  5.4395\\
    5.4410&   5.4408\\
    5.4423&    5.4422\\
    5.4436& 5.4435\\
    5.4449& 5.4448\\
    5.4463&  5.4462\\
    5.4476&  5.4475\\
  5.4489&   5.4475\\
 5.4502&5.4475\\
  5.4515&   5.4501\\
    5.4528& 5.4515\\
 \end{tabular} 
  \begin{tabular} [t] {c c}
  $\Phi_{N}(1)$ & $A N^{-d}$\\   
    5.4542&    5.4541\\
    5.4528    & 5.4554\\
 5.4555&  5.4567\\
   5.4568& 5.4580\\
    5.4581&    5.4593\\
    5.4594&   5.4606\\
    5.4607&5.4620\\
    5.4620&5.4633\\
 5.4633&  5.4646\\
    5.4646& 5.4659\\
   5.4659&  5.4672\\
    5.4672& 5.4685\\
     5.4685&5.4697\\
    5.4698& 5.4710\\
       \end{tabular} 
  \begin{tabular} [t] {c c}
  $\Phi_{N}(1)$ & $A N^{-d}$\\
       5.4711&    5.4723\\
    5.4723& 5.4736\\
    5.4736&  5.4749\\
    5.4749&5.4762\\
  5.4633&  5.4646\\
    5.4646& 5.4659\\
    5.4659&  5.4672\\
    5.4672& 5.4685\\
    5.4685&5.4697\\
    5.4698& 5.4710\\
    5.4711&    5.4723\\
    5.4723& 5.4736\\
    5.4736&  5.4749\\
    5.4749&5.4762\\
    \end{tabular}

      \subsubsection {d=   -0.0250}
  \begin{tabular} [t] {c c}
  $\Phi_{N}(1)$ & $A N^{-d}$\\
    2.3768&  2.3768\\
    2.3769&     2.3769\\
    2.3770&  2.3770\\
    2.3771& 2.3771\\
    2.3772&    2.3772\\
    2.3773&   2.3773\\
    2.3774& 2.3774\\
    2.3775& 2.3775\\
    2.3776&    2.3776\\
    2.3777& 2.3777\\
     2.3778&    2.3778\\
      \end{tabular}
  \begin{tabular} [t] {c c}
  $\Phi_{N}(1)$ & $A N^{-d}$\\
     2.3779&  2.3779\\
    2.3780&  2.3780\\
    2.3781&  2.3781\\
    2.3782& 2.3782\\
    2.3783&    2.3783\\
    2.3784& 2.3784\\
    2.3785& 2.3785\\
    2.3786&    2.3786\\
    2.3787&   2.3787\\
     2.3788& 2.3787\\
    \end{tabular} 
  \begin{tabular} [t] {c c}
  $\Phi_{N}(1)$ & $A N^{-d}$\\
     2.3789&  2.3788\\
    2.3789&    2.3789\\
    2.3790&   2.3790\\
    2.3791&  2.3791\\
    2.3792&2.3792\\
    2.3793& 2.3793\\
    2.3794&   2.3794\\
    2.3795& 2.3795\\
  2.3796& 2.3796\\
   2.3797&   2.3797\\
            \end{tabular}
  \begin{tabular} [t] {c c}
  $\Phi_{N}(1)$ & $A N^{-d}$\\
       2.3798&  2.3798\\
    2.3799&2.3799\\
    2.3800&2.3800\\
    2.3801&  2.3801\\
    2.3802&  2.3802\\
    2.3803&   2.3803\\
    2.3804&   2.3804\\
    2.3805&   2.3805\\
    2.3805& 2.3805\\
    2.3806&    2.3806\\
         \end{tabular}

        \subsubsection {d= 0.100}
  \begin{tabular} [t] {c c}
  $\Phi_{N}(1)$ & $A N^{-d}$\\
    0.9706& 0.9707\\
    0.9704& 0.9705\\
    0.9702& 0.9703\\
   0.9701&  0.9702\\
    0.9699&0.9700\\
    0.9698&0.9699\\
    0.9696&0.9697\\
    0.9695&0.9695\\
    0.9693&0.9694\\
    0.9691&0.9692\\
     0.9690&0.9691\\
      \end{tabular}
        \begin{tabular} [t] {c c}
  $\Phi_{N}(1)$ & $A N^{-d}$\\
         0.9688&0.9689\\
    0.9687&0.9687\\
    0.9685&0.9686\\
    0.9684&0.9684\\
    0.9682&0.9683\\
    0.9680&0.9681\\
    0.9679&0.9679\\
    0.9677&0.9678\\
           0.9676&0.9676\\
     0.9674&0.9675\\
      0.9673&0.9673\\
       \end{tabular} 
  \begin{tabular} [t] {c c}
  $\Phi_{N}(1)$ & $A N^{-d}$\\
       0.9671&0.9672\\
    0.9670&0.9670\\
    0.9668&0.9669\\
    0.9667&0.9667\\
    0.9665&0.9665\\
    0.9664&0.9664\\
      0.9662&0.9662\\
    0.9661&0.9661\\
        0.9659& 0.9659\\
    0.9658&0.9658\\
          0.9656&0.9656\\
     \end{tabular}
      \begin{tabular} [t] {c c}
  $\Phi_{N}(1)$ & $A N^{-d}$\\
          0.9654&0.9655\\
    0.9653& 0.9653\\ 
    0.9651& 0.9652\\
    0.9650&0.9650\\
    0.9648&0.9649\\
    0.9647&0.9647\\
    0.9646&0.9646\\
    0.9644&0.9644\\
   
   \end{tabular} 
     \subsubsection {d= 0.2250}
  \begin{tabular} [t] {c c}
  $\Phi_{N}(1)$ & $A N^{-d}$\\
   0.3566&0.3568\\
    0.3565&0.3567\\
    0.3563&0.3565\\
    0.3562&0.3564\\
    0.3561&0.3563\\
    0.3559&0.3561\\
    0.3558&0.3560\\
    0.3557&0.3559\\
    0.3556&0.3557\\
   0.3554&0.3556\\
    0.3553&0.3555\\
                \end{tabular} 
         \begin{tabular} [t] {c c}
  $\Phi_{N}(1)$ & $A N^{-d}$\\
       0.3552&0.3553\\
    0.3551&0.3552\\
    0.3549&0.3551\\
    0.3548&0.3550\\
    0.3547&0.3548\\
    0.3546&0.3547\\
    0.3544&0.3546\\
    0.3543&0.3544\\
    0.3542&0.3543\\
    0.3541&0.3542\\
        0.3539&0.3540\\
   \end{tabular}
         \begin{tabular} [t] {c c}
  $\Phi_{N}(1)$ & $A N^{-d}$\\
    0.3538&0.3539\\
    0.3537&0.3538\\
    0.3536&0.3537\\
    0.3535&0.3535\\
    0.3533&0.3534\\
    0.3532&0.3533\\
    0.3531&0.3532\\
    0.3530&0.3530\\
    0.3528& 0.3529\\
   0.3527& 0.3528\\
    0.3526& 0.3526\\
   \end{tabular}
         \begin{tabular} [t] {c c}
  $\Phi_{N}(1)$ & $A N^{-d}$\\
      0.3525& 0.3525\\
    0.3524& 0.3524\\
    0.3522& 0.3523\\
    0.3521& 0.3521\\
    0.3520& 0.3520\\
    0.3519& 0.3519\\
    0.3518& 0.3518\\
    0.3517& 0.3517\\
\end{tabular}
 \bibliography{Toeplitzdeux}
\end{document}